%% file: FurtherHopping.tex
\documentclass{amsart}
\usepackage{graphicx}
\usepackage{color}
\newtheorem{thm}{Theorem}[section]

\newtheorem{lem}{Lemma}[section]

\theoremstyle{definition}

\theoremstyle{remark}

\numberwithin{equation}{section}

\newcommand{\regT}{{\rm{T}}}
\newcommand{\regF}{{\rm{F}}}
\newcommand{\regL}{{\rm{L}}}
\newcommand{\regR}{{\rm{R}}}
\begin{document}

\title[Further Hopping with Toads and Frogs]{Further Hopping with Toads and Frogs}%
\author{Thotsaporn ``Aek'' Thanatipanonda}%
\address{Rutgers, The State University of New Jersey}
\address{110 Frelinghuysen Rd Piscataway}
\address{NJ 08854-8019, USA}
\email{thot@math.rutgers.edu}%


\begin{abstract}
    We show the value of positions of the combinatorial game ``Toads and Frogs''.
 We present new values of starting positions. Moreover, we discuss the values 
of all positions with exactly one $\Box, \regT^{a}\Box\Box \regF^{a}, 
\regT^{a} \Box \Box \Box \regF \regF \regF,\regT^{a}\Box\Box \regF^{b}$, 
$\regT^{a}\Box\Box\Box \regF^{b}$. 
At the end, we post five new conjectures and discuss the possible future work.
\end{abstract}
\maketitle
\section{Introduction}

The game {\it Toads and Frogs}, invented by Richard Guy,
is extensively discussed in ``Winning Ways''[1], the famous
classic by Elwyn Berlekemp, John Conway, and Richard Guy, that is
the {\it bible} of combinatorial game theory. \\

This game got so much coverage because of the
simplicity and elegance of its rules, the beauty of its
analysis, and  as an example of a  combinatorial game
whose positions do not always have values that are numbers. \\

The game is played on a $1 \times n$ strip with either
Toad(T) , Frog(F)
or $\Box$ on the squares. Left plays T and Right plays F.
T may move to the immediate square on its right, if it happens
to be empty,
and F moves to the next empty square on the
left, if it is empty.
If T and F are next to each other, they have an option to jump over one another,
in their designated directions, provided they land on an empty square.
(see [1] page 14).\\

\noindent In symbols: the following moves are legal for Toad: \\

\noindent $\dots \regT \Box \dots \;\ \rightarrow \;\ \dots \Box \regT \dots $,\\
$\dots \regT \regF \Box \dots \;\ \rightarrow \;\ \dots \Box \regF \regT \dots $
\;\ , \\

\noindent and  the following moves are legal for Frog: \\

\noindent $\dots \Box \regF \dots \;\ \rightarrow \;\ \dots \regF \Box \dots $,\\
$\dots \Box \regT \regF \dots \;\ \rightarrow \;\ \dots \regF \regT \Box \dots
$. \\

Already in ``Winning Ways''[1] there is
some analysis of Toads and Frogs positions, but on {\it specific}, small
boards, such as $\regT \regT \regT \Box \regF \regF$. In 1996, Erickson[2]
analyzed more general positions. At the end he made five
conjectures about the values of some
families of positions. All of them are starting
positions (positions where all T are rightmost and all F are leftmost).\\

In [3], we discussed the symbolic finite-state approach to prove the value of the positions in class A and B. However,the patterns of the value of Toads and Frogs game are not limited to only class A or class B. There are also patterns in the positions where the variables are on both T and F, for example $\regT^{a} \Box\Box \regF^{b}$. In this paper, we will analyze some of these positions.\\

To be able to understand the present article, reader need a minimum knowledge of combinatorial game theory, that can be found in [1]. In particular, readers should be familiar with the notions of {\it value} of a game and sum games.\\

Let's recall the bypass reversible move rule, dominated options rule and inequality of two games. (see [1] page 33, 62-64).\\

\noindent \textbf{The Bypassing right's reversible move rule}\\

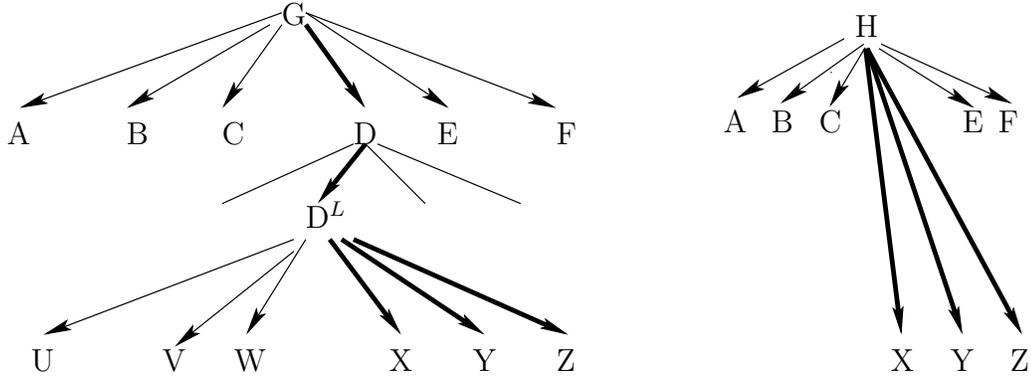
\begin{figure} [ht]
\input{tf10.pstex_t}
\caption{The Bypassing reversible move rule.}
\end{figure}

\noindent $G = H$ if $D^{L} \geq G.$\\

\noindent \textbf{The Dominated options rule} \\

\noindent Let $G = \{A,B,C,... \mid D,E,F,.. \}.$\\
If $A \geq B$ and $D \geq E$ then $G = \{A,C,... \mid E,F,.. \}.$\\

\noindent \textbf{Sum and Inequality of two games}\\

\noindent Let $G = \{ G^{L} \mid G^{R}\}.$\\
Let $H = \{ H^{L} \mid H^{R}\}.$\\

\noindent $G + H = \{ G^{L}+H, G+H^{L} \mid G^{R}+H, G+H^{R} \}.$\\

\noindent $G \geq 0$ if Right goes first and Left wins.\\
$G \leq 0$ if Left goes first and Right wins.\\
$G \geq H$ is equivalent to $G+(-H) \geq 0.$\\

The only notations we use are $*  ( = \{ 0 \mid 0 \})$ and $n* (= \{ n \mid n \})$. We will not use any shorthand notation like $\uparrow$, $\Uparrow$, etc.\\

\section {The general classes A and B}

\noindent \textbf{Definition}:\\
General class Ai: all positions with (numeric) i number of $\Box$ (with symbols on both T and F).\\
General class Bi: all positions with (numeric) i number of F (with symbols on both T and $\Box$). \\

The general classes A and B are the more general positions of the classes A and B that we talked about in [3].\\

For the general class, we can not apply the finite state method that we use in [3] anymore since we now have infinitely many positions that come from the combination of the two letters with symbols on them. However we manage to categorize all positions in the general class A1, the class of all positions with exactly on $\Box$. It is in fact the only general class that we manage to do it for.\\

Many positions in these classes do not have a nice compact formula; for example in A2, $\regT^{a}\Box \regT \regF \Box \regT \regF^{b}$. On the other hand, many positions have a nice formula. We will prove some of the starting positions like $\regT^{a} \Box\Box \regF^{a}, \regT^{a}\Box\Box \regF^{b}$ later on in the appendix of this paper.\\

Once we detect the patterns of the positions, the proof is quite routine. We now do the proof of each specific position by hand with the help of a computer. We hope to see the computer playing more roles in assisting with the proofs in the future.\\
\section {Table}

We present the values of some starting positions in this section. We have a fast program written in Java to calculate the outcome of the sum of two given positions (=,$>$,$<$,$||$). This program does not calculate the value of the sum of two games. It only gives the outcome. It works well with the positions that have a simple value. The author's brother and the author wrote this program originally to check the value of the game of the form $\regT^{a} \Box^{b} \regF^{a}$ for which so far the values of the game are 0 or * except the column b=2 which will be proved to be infinitesimal when $a \geq 4$. The program is on the author's web site and we present the tables here. \\ \\

\newpage

\begin{figure} [htbp]
\begin{tabular}{|c|c|c|c|c|c|c|c|c|c|c|c|c|c|c|c|c|c|c|c|c|}

  \hline
  $a \backslash b$  & 1 & 2 & 3 & 4 & 5 & 6 & 7 & 8 & 9 & 10 & 11 & 12 & 13 & 14 & 15 & 16 & 17 & 18 & 19 & 20 \\
    \hline
  1  & $*$ & $0$ & $*$ & $0$ & $*$ & $0$ & $*$ & $0$ & $*$ & $0$ & $*$ & $0$ & $*$ & $0$ & $*$ & $0$ & $*$ & $0$ & $*$ & $0$ \\
  2  & $*$ & $*$ & $*$ & $*$ & $0$ & $0$ & $*$ & $0$ & $0$ & $0$ & $0$ & $0$ & $*$ & $0$ & $0$ & $0$ & $0$ & $0$ & $*$ & $0$ \\
  3  & $*$ & $\pm \frac{1}{8}$ & $0$ & $*$ & $0$ & $*$ & $0$ & $0$ & $0$ & $0$ & $0$ & $0$ & $0$ & $*$ & $0$ & $0$ & $0$ & $0$ & $0$ & $0$ \\
  4  & $*$ & N & $*$ & $0$ & $0$ & $0$ & $0$ & $*$ & $*$ & $0$ & $0$ & $0$ & $0$ & $0$ & $0$ & $0$ & $0$ & $0$ & $0$ & $0$ \\
  5  & $*$ & N & $*$ & $*$ & $*$ & $0$ & $*$ & $0$ & $0$ & $0$ & $ \neq 0$ & $0$  &&&&&&&&\\
  6  & $*$ & N & $*$ & $*$ & $*$ & $*$ & $*$ & $\neq 0$ &&&&&&&&&&&&\\
  7  & $*$ & N & $*$ & $*$ & $*$ & $\neq 0$ & $\neq 0$  &&&&&&&&&&&&&\\
  8  & $*$ & N & $*$ & $*$ & $*$  &&&&&&&&&&&&&&&\\
  9  & $*$ & N & $*$ & $*$  &&&&&&&&&&&&&&&&\\
  10 & $*$ & N & $*$ & $*$  &&&&&&&&&&&&&&&&\\

  \hline
\end{tabular}
\caption{ $\regT^{a} \Box^{b} \regF^{a}$}
\end{figure}

\noindent Note 1) For b where $21 \leq b \leq 103 ,\regT^{2} \Box^{b} \regF^{2}$ = 0 except b = 25,31,37,43,49,55,61,67,73,79,85,91,97,103.\\
Note 2) For b where $21 \leq b \leq 53 , \regT^{3} \Box^{b} \regF^{3}$ = 0 except b = 29. \\
Note 3) N is an infinitesimal, it is long. We are not writing it out here.\\

\begin{figure} [htbp]
\begin{tabular}{|c|c|c|c|c|c|c|c|c|c|c|c|c|c|c|c|c|}
  \hline
  $a \backslash b$   & 3 & 4 & 5 & 6  \\
    \hline
  1   & $2*$ & $3$ & $4*$ & $5$   \\
  2   & $1*$ & $2*$ & $3$ & $\frac{7}{2}$ \\
  3   & $\{ 1 \mid \{ \frac{1}{2} \mid 0 \} \}$ & $1*$ & $2$ & $\frac{11}{4}$ \\
  4   & $1*$ & $ \{\{ 2* \mid 1*\} \mid\mid \{ \frac{1}{2} \mid 0 \}\} $ & $1$ & $2$ \\
  5   & $1*$ & $1*$  & $ \frac{5}{4} < V < 2$ & $1$ \\
  6   & $1*$ & $1*$  & $2*$ & $1 < V < 2$ \\
  7   & $1*$ & $1*$  & $ || \;\ 2$ & \\
  8   & $1*$ & $1*$  & $ || \;\ 2 $&\\

  \hline
\end{tabular}

\caption{$\regT^{a+1} \Box^{b} \regF^{a}, \mbox{first part}$}
\end{figure}

\begin{figure}
\begin{tabular}{|c|c|c|c|c|c|c|c|c|c|c|c|c|c|c|c|c|}
  \hline
  $a \backslash b$  & 7 & 8 & 9 & 10 & 11 & 12 & 13 & 14 & 15 & 16 \\
    \hline
  1  & $6*$ & $7$ & $8*$ & $9$ & $10*$ & $11$ & $12*$ & $13$ & $14*$ & $15$  \\
  2  & $5*$ & $\frac{11}{2}$ & $\frac{13}{2}$ & $\frac{15}{2}$ & $\frac{17}{2}$ & $\frac{19}{2}$ & $\frac{41}{4}$ & $\frac{23}{2}$ & $12*$ & $13$  \\
  3  & $\frac{15}{4}$ & $\frac{9}{2}$ & $< \frac{11}{2}$ & $< 6$ &  &  &  &  &  & \\
  4  & $ \frac{5}{2} < V < 3$ & $< 4$ &  &  &  &  &  &  &  &   \\
  5  & $ \frac{5}{2} < V < 3$ &  &  &  &  & &&& &  \\
  6  & $ < 2 \mbox{ and } || \;\ 3$ &&&&&&&&&\\
  7  &&&&&&&&&& \\

  \hline
\end{tabular}
\caption {$\regT^{a+1} \Box^{b} \regF^{a}, \mbox{second part}$}
\end{figure}
\newpage

\begin{figure} [htbp]
\begin{tabular}{|c|c|c|c|c|c|c|c|c|c|c|c|c|}
  \hline
  $a \backslash b$   & 3 & 4 & 5 & 6 & 7 & 8 & 9 & 10 & 11 & 12  \\
    \hline
  1   & $4*$ & $6$ & $8*$ & $10$ & $12*$ & $14$ & $16*$ & $18$ & $20*$ & $22$ \\
  2   & $2*$ & $4*$ & $\frac{95}{16}$ & $\frac{15}{2}$ & $9*$ & $\{ 11 \mid 11* \}$ & $12$ & $\frac{29}{2}$ & $15*$ & $17*$ \\
  3   & $\frac{3}{2}$ & $\{ \{ \frac{5}{2} \mid 2 \} \mid\mid 2 \}$ & $4$ & $\{ \frac{11}{2} \mid \frac{41}{8} \}$ & $7$ &&&&& \\
  4   & $2*$ & $\{ \{ 4* \mid 2* \} \mid\mid \{ \frac{3}{2} \mid 1 \} \}$ & $3$ & $4$ & $5*$ &  &  & &  &  \\
  5   & $2*$ & $2*$ & $3<V<4$ & $3$ & $5$ &  &  &  &  &   \\
  6   & $2*$ & $2*$ & $ || \;\ 4 $ &&&&&&&\\
  7   & $2*$ & $2*$  & $<3$ &  &  &&&&& \\

  \hline
\end{tabular}
\caption{$\regT^{a+2} \Box^{b} \regF^{a}$}
\end{figure}

\begin{figure} [htbp]
\begin{tabular}{|c|c|c|c|c|c|c|c|c|c|}
  \hline
  $a \backslash b$  & 3 & 4 & 5 & 6 & 7 & 8    \\
    \hline
  1  & $ 6* $ & $9$ & $12* $ & $15$ & $18*$  & $21$  \\
  2  & $ 3* $ & $\{ 6 \mid \frac{11}{2} \})$ & $\{ \frac{17}{2} \mid 8 \})$ & $ 11* $ & $13$ & $\frac{31}{2}$     \\
  3  & $\frac{5}{2}$ & L  & $\frac{41}{8}$ & $ || \;\ 8$ & &    \\
  4  & $ 3* $ & $\frac{5}{2}$ & $5$ & $ || \;\ 5$ & &   \\
  5  & $ 3* $ & $|| \;\ 3 $ & $ 5 < V < 6$ &  &  &   \\
  6  & $ 3*$ & $ 3* $ &  &  & &   \\

  \hline
\end{tabular}
\caption{$\regT^{a+3} \Box^{b} \regF^{a}$}
\bigskip
\flushleft
Note 1) $||G$ means ``can not be compared to G''.\\
Note 2) We drop the values of the first two columns where $b=1,2$ since they will all be proved in the appendix.\\
Note 3) L means long. We are not writing it out here.\\
\end{figure}

\noindent Erickson's conjecture 4 is false since $\regT^{7} \Box^{7} \regF^{6} > 2$.\\

\noindent (\textbf{Erickson's conjecture 4}: $\regT^{a} \Box^{a} \regF^{a-1} = 1$ or $\{ 1 \mid 1 \}$ for all $a \geq 1$.) \\

The author believes there are no patterns in $a$ for positions of the form $\regT^{a+k} \Box^{a+l} \regF^{a}$; for any fixed $k \geq 1, l \geq 0$.\\
\section {New Conjectures and Future Work}

In [2], Jeff Erickson made 6 conjectures. Jesse Hull proved conjecture 6 (Toads and Frogs is NP-hard) in 2000. Doron Zeilberger and the author of this paper proved conjecture 1 (later on in this paper), 2 (in [3]), 3 (later on in this paper) and disprove conjecture 4. Conjecture 5 is still open. We restate conjecture 5 here.\\

\noindent \textbf{Erickson's conjecture 5}:\\
$\regT^{a} \Box^{b} \regF^{a}$ is an infinitesimal for all a,b except (a,b) =(3,2) \\

This conjecture seems very interesting and hard but not impossible to prove. We split Erickson's fifth conjecture
into 2 stronger conjectures which are conjectures 3 and 4 here.\\

 We believe there are still a lot of nice patterns and conjectures in this game that we overlooked. Once RAM gets cheaper and Maple gets faster, we will have more information.\\

\noindent \textbf{Conjecture 1)} Assume $b \geq 0$, $a \geq 1$, $L \geq 0$ and $R \geq 0$ \\

\noindent 1.1)$\Box^{R}\regT^{a}\Box^{b} \regF \Box^{R} =  \left\{ \begin{array}{ll}
                                 \{\{a-2 \mid 1 \} \mid 0 \} & \mbox{if} \;\ R=0 \;\ \mbox{and } \;\ b=1 \\
                                (a-1)(b-1+R)  & \mbox{if} \;\ b \;\ \mbox{is even} \\
                                (a-1)(b-1+R)* & \mbox{if} \;\ b \;\ \mbox{is odd and} (R,b) \neq (0,1)
                                \end{array}
                                \right. $  \\ \\

\noindent 1.2) For $R \geq 1$, $\Box^{R-1}\regT^{a}\Box^{b}\regF\Box^{R} = \left\{ \begin{array}{ll}
                                         (a-1)(b-1+R) & \mbox{if} \;\ b \;\ \mbox{is even} \\
                                         1/2+(a-1)(b-1+R) & \mbox{if} \;\ b \;\ \mbox{is odd} \\
                                          \end{array}
                                          \right. $  \\ \\

\noindent 1.3) For $R-L \geq  2, \Box^{L}\regT^{a}\Box^{b}\regF\Box^{R} = (R-L-1)+(a-1)(b-1+R)$\\

\noindent \textbf{Conjecture 2)} For $a \geq 7$,         $\regT \regT\Box^{a}\regF\regF =  \left\{ \begin{array}{ll}
                              * & \mbox{when } a = 7+6n, \;\ n \geq 0 \\
                              0 & \mbox{otherwise}. \\
                              \end{array}
                              \right. $  \\

\noindent \textbf{Conjecture 3)} $\regT^{a}\Box^{b} \regF^{a} = *$ for any $a> b>0$, except for $b=2$. \\

\noindent \textbf{Conjecture 4)} $\regT^{a}\Box^{b} \regF^{a} = 0$ or $*$ for any $b \geq a > 0 $.\\

\noindent \textbf{Conjecture 5)} For a fixed integer $C \geq 3$, $\exists a_0$ such that $\regT^{C}\Box^{a} \regF^{C} = 0$ for all $a \geq a_0$.\\

\noindent \textbf{Future Work}\\
1) Categorize all the positions that have exactly one F (general class B1) (conjecture 1 might be a good start).\\

\section{Acknowledgements}

I want to thank my twin brother Thotsaphon Thanatipanonda for writing the fast program in Java
to find the counterexample of conjecture 4 of Erickson. I also want to thank my advisor Doron Zeilberger for always
supporting and motivating. There will not be this paper without him.\\
\renewcommand{\theequation}{A-\arabic{equation}}
\setcounter{equation}{0}  

\section*{APPENDIX}  
\appendix{}

\section{Positions with one $\Box$}

In this section we classify all the positions that have one $\Box$. The lemmas below give a recurrence to the positions. We can compute the values of the positions in polynomial time. We omit the proofs here.\\

\noindent \textbf{Notation} \\
$O(x) = \{0 \mid x \}.$ \\
$O^{a}(x) = O(...(O(O(x)))) \;\ a$ times. \\
$\widetilde{\regL}$ = any combination of $\regT$ and $\regF$ that has $\regF$ as its right most entry. For example $\regT \regT \regF \regT \regF$.\\
$\widetilde{\regR}$ = any combination of $\regT$ and $\regF$ that has $\regT$ as its left most entry. For example $\regT \regF \regF \regT \regF$. \\

\begin{lem} Death Leaps Principle(DLP): The position with one empty square for which the only possible move for both sides is a jump has value 0. For example $\regT \regF \regT \regT \regF\Box \regT \regF \regT \regF \regF$.
\end{lem}

\begin{lem} $\widetilde{\regL}\regT \Box \regF \widetilde{\regR} = * .$\end{lem}
\begin{lem} $\widetilde{\regL}\regT^{a} \Box \regF^{b} \widetilde{\regR} = * , \;\ a \geq 2, b \geq 2.$\end{lem}

\begin{lem} $\widetilde{\regL}\regT^{a} \Box (\regT\regF)^{b}\regT\regF^{c}\widetilde{\regR} = \{ a-1 \mid (O^{b}(\widetilde{\regL}\regT^{a}\regF(\regT\regF)^{b}\regT \Box \regF^{c-1}\widetilde{\regR}) \}, \;\ a \geq 1, b\geq 0, c\geq 2.$ \end{lem}

\noindent $\widetilde{\regL}\regT^{a} \Box \regF (\regT\regF)^{b}\regT\regF^{c} \widetilde{\regR} = \{\{a-2 \mid O^{b+1}(\widetilde{\regL}\regT^{a-1}\regF(\regT\regF)^{b+1}\regT \Box \regF^{c-1}\widetilde{\regR})\} \mid \mid 0 \}, \;\ a \geq 2, b\geq 0, c\geq 2.$ \\

\noindent \textbf{Example 1:} $\regT^{a} \Box \regF (\regT\regF)^{b}\regT\regF^{2} =  \{\{a-2 \mid O^{b+1}(*)\} \mid \mid 0 \}, \;\ a \geq 1, b \geq 0.$ \\

\noindent \textbf{Example 2:} $\regT^{3} \Box \regF (\regT\regF)^{b}\regT\regF^{c} =  \{\{1 \mid O^{b+1}(\regT^{2}\regF(\regT\regF)^{b+1}\regT \Box \regF^{c-1} )\} \mid \mid 0 \}, \;\ b \geq 0, c \geq 2.$ \\

\noindent \textbf{Example 3:} $\regT^{a} \Box \regF (\regT\regF)^{b}\regT\regF^{4} =  \{\{a-2 \mid O^{b+1}(\regT^{a-1}\regF(\regT\regF)^{b+1}\regT \Box \regF^{3})\} \mid \mid 0 \}, \;\ a \geq 3, b \geq 0.$ \\

\noindent \textbf{Note)} We get the implicit value of example 2 from example 1 and implicit value of example 3 from example 2. We will get the value recursively this way.\\

\begin{lem}
 $\widetilde{\regL}\regT^{a} \Box (\regT\regF)^{b} = \{ a-1 \mid (\frac{1}{2})^{b-1}\}, \;\ a \geq 1, b\geq 1.$
 \end{lem}

\begin{lem} $\widetilde{\regL}\regT^{a} \Box \regF (\regT\regF)^{b} = \{\{a-2 \mid (\frac{1}{2})^b\} \mid \mid 0 \} , \;\ a \geq 2, b \geq 0.$\end{lem}

\begin{lem} $\widetilde{\regL}\regT^{a} \Box \regF (\regT\regF)^{b}\regT\regF^{c} \widetilde{\regR} = \regT^{a} \Box \regF (\regT\regF)^{b}\regT\regF^{c}, \;\ b \geq 0 $ \end{lem}
\noindent when ($c$ is even and $a \geq c-1$) or ($a$ is odd and $c \geq a-1$).\\

\noindent \textbf{Note:} \\

1) When $a$ is even ,$a \geq 2$ and $c$ is odd ,$c \geq 3$,
The recursive is going to bounce back and forth between positions.\\
We will refer to $\widetilde{\regR}$ if $a>c$. We will refer to $\widetilde{\regL}$ if $c>a$.
Then the positions will start over again.\\

2) When  $a$ is even , $c$ is even and $a < c-1$ then we will refer to $\widetilde{\regL}$.\\
When  $a$ is odd , $c$ is odd and $c < a-1$ then we will refer to 
$\widetilde{\regR}$.\\

\section{ Lemma and Convention.}

\noindent We will refer to the lemma below a lot. We state it here.\\

\noindent \begin{lem} \textbf{One side Death Leap Principle} (One side DLP): if X is the position where the only possible move of Left is a jump and there is no two or more consecutive empty square in X then X $\leq 0$.\end{lem}

\noindent \textbf{Proof:} We have to show that when Left moves first and 
two 
players take turn playing, Left will lost(Left will run out of the legal move first). This is true since after Left jumps over one of the $\regF$, Right can response by moving to the empty square where the $\regF$ was jumped over.  \\

\noindent \textbf{Example 1)}$ \regT \regT  \regF \Box \regT \regT \regF  \Box \regF \leq 0.$\\

\noindent \textbf{Example 2)} $\regT \regT \regT  \regF \Box \regF  \Box \regT \regF \leq 0.$\\

\noindent \textbf{Convention}:\\

In the following sections of the appendix, we will prove the positions using the ``shorthand'' notation. We explain by using an example.\\

\noindent \textbf{Example}: To show: \;\ $ \regT^{a} \regF \Box \regT^{k} \regF \regT \Box \regF^{b} \leq \frac{1}{2} ; \;\ k \geq 0, a \geq 0, b\geq 1.$\\

\noindent We will have to show $\regT^{a} \regF \Box \regT^{k} \regF \regT \Box \regF^{b} - \frac{1}{2} \leq 0.$\\

\noindent That is to show $\regT^{a} \regF \Box \regT^{k} \regF \regT \Box \regF^{b} - \{ 0 \mid 1 \} \leq 0.$\\

To show $G \leq 0$  we need to show that when Left moves first and two players take turn playing, Right will win.
(On the other hand, to show $G \geq 0$  we need to show that when Right moves first and two players take turn playing, Left will win.) \\

We will show that in these two sum games, for all the possible choices of Left move, Right can find a response to the move so that he will win at the end.\\

We will do some case analysis here. In the above position Left has three choices. \\
In the proof we will see \\

$\stackrel {\stackrel{2}{\rightarrow}}{\regT^{a}} \regF \Box \regT^{k} \regF \stackrel {\stackrel{1}{\rightarrow}}{\regT} \Box \regF^{b} \leq \stackrel {\stackrel{3}{\leftarrow}}{\frac{1}{2}}. $\\

\noindent We will write the response of Right immediately. You will see \\

\noindent \textbf{Case 1:} $\regT^{a} \regF \Box \regT^{k} \regF \Box \regT \regF^{b} \leq 0.$ \\

\noindent \textbf{Explanation:} Right responses by picking the left option of $\{ 0 \mid 1 \}$ on the right hand side.\\

\noindent \textbf{Case 2:} $\regT^{a-1} \regF \Box \regT^{k+1} \regF \regT \Box  \regF^{b} \leq \frac{1}{2}.$ \\

\noindent \textbf{Explanation:} Right responses by moving the left most F.\\

\noindent \textbf{Case 3:} $\regT^{a} \regF \Box \regT^{k} \regF \regT \regF \Box \regF^{b-1} \leq 1.$ \\

\noindent \textbf{Explanation:} Left picks the right option of $\{ 0 \mid 1 \}$ on the right hand side. Right responses by moving the right most F.\\

\noindent \textbf{Note:} \\

1) When the position simplify to the known one in [3] or [5] (which is mostly positions in class A), we will claim such result without proving them again.\\

2) The positions we considered in Toads and Frogs are a hot game which means the players select a good move and fight for an advantage. We will not consider the possible move in a cold game which is the whole integer. \\
\section{$\regT^{a} \Box\Box \regF^{a}$}

 We show $\regT^{a} \Box\Box \regF^{a}$ is an infinitesimal , $a \geq 4.$  The observation comes from figure 2 when $b=2$.

\begin{lem}  For any fixed integer $n \geq 3$ , $\widetilde{\regL} \stackrel {\stackrel{2}{\rightarrow}}{\regT^{a}} \regF \Box \regT^{k} \regF \stackrel 
{\stackrel{1}{\rightarrow}}\regT \Box \regF^{b} \leq  \stackrel {\stackrel{3}{\leftarrow}}{\frac{1}{2^{n}}}$, $k \geq 0, a \geq 0 , b \geq 1.$ \end{lem}

\noindent \textbf{Prove:} By induction on $a$. \\

\noindent \textbf{Base Case: $a = 0$}, $\Box \regT^{k} \regF \stackrel {\stackrel{1}{\rightarrow}}\regT \Box \regF^{b} \leq  \stackrel {\stackrel{2}{\leftarrow}}{\frac{1}{2^{n}}}$ . \\

\noindent \textbf{Case 1:} $\Box \regT^{k} \regF \Box \regT \regF^{b} \leq 0$, true by one side DLP.\\

\noindent \textbf{Case 2:} $\Box \regT^{k} \regF \regT \regF \Box \regF^{b-1} \leq \frac{2}{2^{n}}$. The left hand side is $\leq 0$ by one side DLP.\\

\noindent \textbf{Induction Step:} $\widetilde{\regL} \stackrel {\stackrel{2}{\rightarrow}}{\regT^{a}} \regF \Box \regT^{k} \regF \stackrel {\stackrel{1}{\rightarrow}}\regT \Box \regF^{b} \leq  \stackrel {\stackrel{3}{\leftarrow}}{\frac{1}{2^{n}}}$. \\

\noindent \textbf{Case 1:} $\widetilde{\regL}\regT^{a} \regF \Box \regT^{k} \regF \Box \regT \regF^{b} \leq 0$, true by one side DLP.\\

\noindent \textbf{Case 2:} $\widetilde{\regL}\regT^{a-1} \regF \Box \regT^{k+1} \regF \regT \Box \regF^{b} \leq \frac{1}{2^{n}}$, true by induction.\\

\noindent \textbf{Case 3:} $\widetilde{\regL}\regT^{a} \regF \Box \regT^{k} \regF \regT \regF \Box \regF^{b-1} \leq \frac{2}{2^{n}}$. The left hand side is $\leq 0$ by one side DLP.\\

\begin{thm} $\regT^{a} \Box\Box \regF^{a}$ is an infinitesimal , $a \geq 4.$\end{thm}

\noindent \textbf{By symmetry we only need to show:}

\noindent For any fixed integer $n \geq 3, \;\ \regT^{a} \Box\Box \regF^{a}  \leq \frac{1}{2^{n}}, \;\ a \geq 4$.\\

\noindent  $\stackrel {\stackrel{I}{\rightarrow}}{\regT^{a}} \Box\Box \regF^{a}  \leq \stackrel {\stackrel{II}{\leftarrow}}{\frac{1}{2^{n}}}$. \\

\noindent \textbf{I)} $\stackrel {\stackrel{2}{\rightarrow}}{\regT^{a-1}} \Box \stackrel {\stackrel{1}{\rightarrow}}{\regT} \regF \Box \regF^{a-1}  \leq \stackrel {\stackrel{3}{\leftarrow}}{\frac{1}{2^{n}}}$ \\
\textbf{II)} $\stackrel {\stackrel{1}{\rightarrow}}{\regT^{a}} \Box \regF \Box \regF^{a-1}  \leq \stackrel {\stackrel{2}{\leftarrow}}{\frac{2}{2^{n}}}$ \\

\noindent \textbf{I)} {Case 1:} $\stackrel {\stackrel{1}{\rightarrow}}{\regT^{a-1}} \Box \regF \Box \regT \regF^{a-1}  \leq \stackrel {\stackrel{2}{\leftarrow}}{\frac{1}{2^{n}}}$\\

\hspace{3mm} Case 1.1: $\regT^{a-2} \regF \stackrel {\stackrel{1}{\rightarrow}}{\regT} \Box\Box  \regT \regF^{a-1}  \leq \stackrel {\stackrel{2}{\leftarrow}}{\frac{1}{2^{n}}}$ \\

\hspace{8mm} {Case 1.1.1:} $\regT^{a-2} \regF \Box \regT \regF \regT \Box \regF^{a-2}  \leq \frac{1}{2^{n}}$, true by lemma C.1.\\

\hspace{8mm} {Case 1.1.2:} $\regT^{a-2} \regF \stackrel {\stackrel{2}{\rightarrow}}{\regT} \Box \regF \stackrel {\stackrel{1}{\rightarrow}}{\regT} \Box \regF^{a-2}  \leq \stackrel {\stackrel{3}{\leftarrow}}{\frac{2}{2^{n}}}$ \\

\hspace{13mm} {Case 1.1.2.1:} $\regT^{a-2} \regF \stackrel {\stackrel{1}{\rightarrow}}{\regT} \regF \Box\Box \regT \regF^{a-2}  \leq \stackrel {\stackrel{2}{\leftarrow}}{\frac{2}{2^{n}}}$ \\

\hspace{18mm} {Case 1.1.2.1.1:} $\Box \regT \Box \regT \regF^{a-2}  \leq \frac{2}{2^{n}}$.

\hspace{22mm} The left hand side  is $\{ 0 \mid\mid \{ 0 \mid \{ -1 \mid 5-a \} \} \}$. \\

\hspace{18mm} {Case 1.1.2.1.2:} $\regT^{a-2} \regF \regT \regF \Box \regF \regT \Box \regF^{a-3}  \leq \frac{4}{2^{n}}$, true by lemma C.1. \\

\hspace{13mm} {Case 1.1.2.2:} $\regT^{a-2} \regF\Box \regT \regF \regT \regF \Box \regF^{a-3}  \leq \frac{2}{2^{n}}$.

\hspace{17mm} The left hand side  is $\leq 0$ by one side DLP.\\

\hspace{13mm} {Case 1.1.2.3:} $\regT^{a-2} \regF \stackrel {\stackrel{2}{\rightarrow}}{\regT} \Box \regF \stackrel {\stackrel{1}{\rightarrow}}{\regT} \regF \Box \regF^{a-3}  \leq \stackrel {\stackrel{3}{\leftarrow}}{\frac{4}{2^{n}}}$ \\

\hspace{18mm} {Case 1.1.2.3.1:} $\regT^{a-2} \regF \regT \Box \regF \Box \regF \regT \regF^{a-3}  \leq 0$ \\

\hspace{38mm} $\Rightarrow \regT^{a-2} \regF \Box \regT \regF \regF \Box \regT \regF^{a-3}  \leq 0$

\hspace{42mm}The statement is true by one side DLP.\\

\hspace{18mm} {Case 1.1.2.3.2:} $\regT^{a-2} \regF \Box \regT \regF \regT \regF \regF \Box \regF^{a-4}  \leq \frac{4}{2^{n}}$.

\hspace{22mm} The left hand side  is $\leq 0$ by one side DLP. \\

\hspace{18mm} {Case 1.1.2.3.3:} $\regT^{a-2} \regF \regT \regF \Box \regT \regF \Box \regF^{a-3}  \leq \frac{8}{2^{n}}$.

\hspace{22mm} The left hand side  is $\leq 0$ by one side DLP. \\

\hspace{3mm} {Case 1.2:} $\stackrel {\stackrel{2}{\rightarrow}}{\regT^{a-1}} \Box \regF \regF \stackrel {\stackrel{1}{\rightarrow}}{\regT}\Box \regF^{a-2}  \leq \stackrel {\stackrel{3}{\leftarrow}}{\frac{2}{2^{n}}}$ \\

\hspace{8mm} {Case 1.2.1:} $\regT^{a-1} \regF \Box \regF \Box \regT \regF^{a-2}  \leq \frac{2}{2^{n}}$.

\hspace{12mm} The left hand side  is $\leq 0$ by one side DLP.\\

\hspace{8mm} {Case 1.2.2:} $\regT^{a-2} \regF \regT \Box \regF \regT \Box \regF^{a-2}  \leq \frac{2}{2^{n}}$. This is the case 1.1.2 \\

\hspace{8mm} {Case 1.2.3:} $\regT^{a-1} \regF \Box \regF \regT \Box \regF^{a-2}  \leq \frac{4}{2^{n}}$, true by lemma C.1.\\

\noindent {Case 2:} $\regT^{a-2} \Box \regT \regT \regF \regF \Box \regF^{a-2}  \leq \frac{1}{2^{n}}$. The left hand side is 0.\\

\noindent {Case 3:} $\stackrel {\stackrel{1}{\rightarrow}}{\regT^{a-1}} \Box \regT \regF \regF \Box \regF^{a-2}  \leq \stackrel {\stackrel{2}{\leftarrow}}{\frac{2}{2^{n}}}$ \\

\hspace{3mm} {Case 3.1:} $ \regT^{a-2} \Box \regT \regT \regF \regF \Box \regF^{a-2} \leq 0. $ This is clearly true.\\

\hspace{3mm} {Case 3.2:} $\regT^{a-1} \regF \regT \Box \regF \Box \regF^{a-2}  \leq \frac{4}{2^{n}}$. The left hand side  is $\leq 0$.\\

\noindent \textbf{II)}  {Case 1:} $\regT^{a-1} \Box \regT \regF \regF \Box \regF^{a-2}  \leq \frac{2}{2^{n}}$. This is I) case 3.\\

\noindent {Case2:} $\stackrel {\stackrel{1}{\rightarrow}}{\regT^{a}} \Box \regF \regF \Box \regF^{a-2}  \leq \stackrel {\stackrel{2}{\leftarrow}}{\frac{4}{2^{n}}}$\\

\hspace{3mm} {Case 2.1:} $\regT^{a-1} \regF \regT \Box \regF \Box \regF^{a-2}  \leq \frac{4}{2^{n}}$. This is I) case 3.2.\\

\hspace{3mm} {Case 2.2:} $\regT^{a} \regF \Box \regF \Box \regF^{a-2}  \leq \frac{8}{2^{n}}$. The left hand side  is $\leq 0$ by one side DLP.\\

\noindent The theorem is proved. $\Box$ \\
\section{ Three More Results}
\begin{thm}$\regT^{a} \Box\Box\Box \regF\regF\regF = a - \frac{7}{2}$ , $a \geq 5.$
\end{thm}

\begin{thm}  $\regT^{a} \Box \Box \regF^{b} = \{\{a-3 \mid a-b \}\mid \{*\mid 3-b\}\},\;\ a > b \geq 2.$
\end{thm}

\begin{thm} $\regT^{a} \Box\Box\Box \regF^{b} = \{a-b \mid a-b \}$ , \;\ $a \geq 4, b \geq 4$ .\end{thm}

\noindent The proofs are quite tedious. We will show them in the supplement paper [4]. \\

\end{document}

%% file: tf10.pstex_t
\begin{picture}(0,0)%
\includegraphics{tf10.pstex}%
\end{picture}%
\setlength{\unitlength}{3947sp}%
\begingroup\makeatletter\ifx\SetFigFont\undefined%
\gdef\SetFigFont#1#2#3#4#5{%
  \reset@font\fontsize{#1}{#2pt}%
  \fontfamily{#3}\fontseries{#4}\fontshape{#5}%
  \selectfont}%
\fi\endgroup%
\begin{picture}(6423,2356)(4576,-3860)
\put(6451,-2911){\makebox(0,0)[lb]{\smash{{\SetFigFont{12}{14.4}{\rmdefault}{\mddefault}{\updefault}{\color[rgb]{0,0,0}${\rm D}^{L}$}%
}}}}
\put(6301,-1636){\makebox(0,0)[lb]{\smash{{\SetFigFont{12}{14.4}{\rmdefault}{\mddefault}{\updefault}{\color[rgb]{0,0,0}$ \rm G$}%
}}}}
\put(4576,-2386){\makebox(0,0)[lb]{\smash{{\SetFigFont{12}{14.4}{\rmdefault}{\mddefault}{\updefault}{\color[rgb]{0,0,0}$ \rm A$}%
}}}}
\put(5326,-2386){\makebox(0,0)[lb]{\smash{{\SetFigFont{12}{14.4}{\rmdefault}{\mddefault}{\updefault}{\color[rgb]{0,0,0}$ \rm B$}%
}}}}
\put(5926,-2386){\makebox(0,0)[lb]{\smash{{\SetFigFont{12}{14.4}{\rmdefault}{\mddefault}{\updefault}{\color[rgb]{0,0,0}$ \rm C$}%
}}}}
\put(6751,-2386){\makebox(0,0)[lb]{\smash{{\SetFigFont{12}{14.4}{\rmdefault}{\mddefault}{\updefault}{\color[rgb]{0,0,0}$ \rm $D}%
}}}}
\put(7276,-2386){\makebox(0,0)[lb]{\smash{{\SetFigFont{12}{14.4}{\rmdefault}{\mddefault}{\updefault}{\color[rgb]{0,0,0}$ \rm E$}%
}}}}
\put(8026,-2386){\makebox(0,0)[lb]{\smash{{\SetFigFont{12}{14.4}{\rmdefault}{\mddefault}{\updefault}{\color[rgb]{0,0,0}$ \rm F$}%
}}}}
\put(4726,-3811){\makebox(0,0)[lb]{\smash{{\SetFigFont{12}{14.4}{\rmdefault}{\mddefault}{\updefault}{\color[rgb]{0,0,0}$ \rm U $}%
}}}}
\put(5551,-3811){\makebox(0,0)[lb]{\smash{{\SetFigFont{12}{14.4}{\rmdefault}{\mddefault}{\updefault}{\color[rgb]{0,0,0}$ \rm V $}%
}}}}
\put(6001,-3811){\makebox(0,0)[lb]{\smash{{\SetFigFont{12}{14.4}{\rmdefault}{\mddefault}{\updefault}{\color[rgb]{0,0,0}$ \rm W$}%
}}}}
\put(6976,-3811){\makebox(0,0)[lb]{\smash{{\SetFigFont{12}{14.4}{\rmdefault}{\mddefault}{\updefault}{\color[rgb]{0,0,0}$\rm X$}%
}}}}
\put(7501,-3811){\makebox(0,0)[lb]{\smash{{\SetFigFont{12}{14.4}{\rmdefault}{\mddefault}{\updefault}{\color[rgb]{0,0,0}$\rm Y$}%
}}}}
\put(8026,-3811){\makebox(0,0)[lb]{\smash{{\SetFigFont{12}{14.4}{\rmdefault}{\mddefault}{\updefault}{\color[rgb]{0,0,0}$\rm Z$}%
}}}}
\put(9901,-1711){\makebox(0,0)[lb]{\smash{{\SetFigFont{12}{14.4}{\rmdefault}{\mddefault}{\updefault}{\color[rgb]{0,0,0}H}%
}}}}
\put(9076,-2311){\makebox(0,0)[lb]{\smash{{\SetFigFont{12}{14.4}{\rmdefault}{\mddefault}{\updefault}{\color[rgb]{0,0,0}A}%
}}}}
\put(9376,-2311){\makebox(0,0)[lb]{\smash{{\SetFigFont{12}{14.4}{\rmdefault}{\mddefault}{\updefault}{\color[rgb]{0,0,0}B}%
}}}}
\put(9676,-2311){\makebox(0,0)[lb]{\smash{{\SetFigFont{12}{14.4}{\rmdefault}{\mddefault}{\updefault}{\color[rgb]{0,0,0}C}%
}}}}
\put(10576,-2311){\makebox(0,0)[lb]{\smash{{\SetFigFont{12}{14.4}{\rmdefault}{\mddefault}{\updefault}{\color[rgb]{0,0,0}E}%
}}}}
\put(10801,-2311){\makebox(0,0)[lb]{\smash{{\SetFigFont{12}{14.4}{\rmdefault}{\mddefault}{\updefault}{\color[rgb]{0,0,0}F}%
}}}}
\put(10126,-3811){\makebox(0,0)[lb]{\smash{{\SetFigFont{12}{14.4}{\rmdefault}{\mddefault}{\updefault}{\color[rgb]{0,0,0}X}%
}}}}
\put(10501,-3811){\makebox(0,0)[lb]{\smash{{\SetFigFont{12}{14.4}{\rmdefault}{\mddefault}{\updefault}{\color[rgb]{0,0,0}Y}%
}}}}
\put(10876,-3811){\makebox(0,0)[lb]{\smash{{\SetFigFont{12}{14.4}{\rmdefault}{\mddefault}{\updefault}{\color[rgb]{0,0,0}Z}%
}}}}
\end{picture}%